\documentclass{elsart3-1}

\usepackage[english, frenchb]{babel} 

\usepackage[applemac]{inputenc}
\usepackage{amsfonts,amsmath}
\usepackage{latexsym}
\usepackage{amssymb}
\usepackage{euscript}
\usepackage{graphicx}

\newtheorem{theoreme}{Th\'eor\`eme}[section]

\renewcommand{\theequation}{\arabic{equation}}
\def\og{\leavevmode\raise.3ex\hbox{$\scriptscriptstyle\langle\!\langle$~}}
\def\fg{\leavevmode\raise.3ex\hbox{~$\!\scriptscriptstyle\,\rangle\!\rangle$}}

\newcommand{\rz}[1]{\mathbb{#1}}
\def\R{\rz R}  
\def\N{\rz N}  

\def\indi{\mbox{1\hspace{-.25em}I}}

\begin{document}

\begin{frontmatter}
\selectlanguage{francais}
\title{Quelques approximations du temps local brownien}
\author[iecn]{Blandine B\'erard Bergery},
\ead{berardb@iecn.u-nancy.fr}
\author[iecn]{Pierre Vallois}
\ead{Pierre.Vallois@iecn.u-nancy.fr}
\address[iecn]{Universit\'e Henri Poincar\'e, Institut de
Math\'ematiques Elie Cartan, B.P. 239, F-54506 Vand\oe uvre-l\`es-Nancy Cedex, France}

\begin{abstract}
\selectlanguage{frenchb}
On définit plusieurs approximations du
processus des temps locaux $(L_t^x)_{t\geqslant 0}$ au niveau $x$
du mouvement brownien réel $(X_t)$. En particulier, on montre que
$ \frac{2}{\epsilon}\int_0^{t} X_{(u+\epsilon)\wedge t}^+ \indi_{ \{X_u \leqslant 0\} } du +
\frac{2}{\epsilon}\int_0^{t} X_{(u+\epsilon) \wedge t}^- \indi_{ \{X_u>0\} } du$ et $\frac{4}{\epsilon}\int_0^{t} X_u^- \indi_{ \{X_{(u+\epsilon) \wedge  t} > 0\} } du$ convergent au sens ucp vers $L_t^0$,
  lorsque $\epsilon \to 0$. D'autre part, on montre que $ \frac{1}{\epsilon}\int_0^t ( \indi_{\{ x<X_{s+\epsilon}\}} - \indi_{\{ x<X_{s}\}} )( X_{s+\epsilon}-X_{s} )ds$ converge vers $L_t^x$  dans $L^2(\Omega)$ et que la vitesse de convergence est d'ordre $\epsilon^\alpha$, pour tout $\alpha < \frac{1}{4}$. 
\vskip 0.5\baselineskip
\selectlanguage{english}
\noindent{\bf Abstract}
\vskip 0.5\baselineskip
\noindent
{\bf Some Brownian local time approximations} \\
\noindent We give some approximations of  the local time process $(L_t^x)_{t\geqslant 0}$ at level $x$ of the real Brownian motion $(X_t)$.
 We prove  that $ \frac{2}{\epsilon}\int_0^{t} X_{(u+\epsilon)\wedge t}^+ \indi_{ \{X_u \leqslant 0\} } du +
\frac{2}{\epsilon}\int_0^{t} X_{(u+\epsilon) \wedge t}^- \indi_{ \{X_u>0\} } du$ 
and $\frac{4}{\epsilon}\int_0^{t} X_u^- \indi_{ \{X_{(u+\epsilon) \wedge  t} > 0\} } du$ converge  in the ucp sense to $L_t^0$, as $\epsilon \to 0$. 
 We show  that $ \frac{1}{\epsilon}\int_0^t ( \indi_{\{ x<X_{s+\epsilon}\}} - \indi_{\{ x<X_{s}\}} )   ( X_{s+\epsilon}-X_{s} )ds$ 
 goes  to $L_t^x$  in $L^2(\Omega)$   as $\epsilon \to 0$, and that the rate of convergence is of order $\epsilon^\alpha$, for any $\alpha < \frac{1}{4}$. 

\vskip 0.5\baselineskip
\noindent{\bf Mots-clés}: temps local, intégration stochastique
 par régularisation, variation quadratique, vitesse de convergence.\\
\noindent{\bf classification AMS}:  60G44, 60H05, 60H99, 60J55, 60J65.
\end{abstract}
\end{frontmatter}

\selectlanguage{frenchb} 
\noindent Dans cette Note, le processus $X$ est continu, et la convergence en probabilité,
uniformément sur les compacts, sera notée ucp (voir Section II.4 de \cite{16}).

\section{D\'efinition  du premier schéma d'approximation} \label{sec1}
\noindent\textbf{1.1.}  Il existe déjà de nombreuses approximations du temps local de larges classes de processus (voir  \cite{13}, \cite{17}, \cite{1}). 
L'objectif de cette Note est de présenter de nouveaux schémas d'approximation du temps local du
mouvement brownien standard réel et de certaines martingales browniennes. On se place dans le cadre de l'intégration par régularisation définie par Russo et Vallois ( \cite{4}, \cite{6}, \cite{8}). On rappelle (c.f.  \cite{6})  que la covariation $[X,Y]$ est la limite au sens ucp de $ \frac{1}{\epsilon}\int_0^t
 \left( Y_{s+\epsilon} - Y_s \right) \left( X_{s+\epsilon}-X_{s}  \right)ds$,  si cette limite existe. On définit
  \begin{equation}\label{definitionJ}
 J_{\epsilon}(t,y)=  \frac{1}{\epsilon}\int_0^t
 \left( \indi_{\{ y<X_{s+\epsilon}\}} - \indi_{\{ y<X_{s}\}} \right)
   \left( X_{s+\epsilon}-X_{s}  \right)ds, \quad  y \in \mathbb{R}, t \geqslant 0.
\end{equation}

\noindent \textbf{1.2.} Il est facile de montrer que  si $f\in\mathcal C^0(\R)$ et $X$ admet une variation quadratique $[X,X]$, alors
\begin{equation}
\label{consequenceJ}
 \lim_{\epsilon \to 0}\textrm{ (ucp) } \int_{\R} f(y)J_{\epsilon}(t,y) dy = \int_0^t f(X_s)d[X,X]_s, \quad f\in\mathcal C^0(\R).
\end{equation}
Lorsque $X$ est  une semi-martingale, $[X,X]$ est égal à  la variation quadratique usuelle et $X$ a un processus des temps locaux  $(L_t^a)_{a\in \R, t \geqslant 0}$. La formule de densité d'occupation permet d'écrire (\ref{consequenceJ}) sous la forme:
\begin{equation*}
 \lim_{\epsilon \to 0}\textrm{ (ucp) } \int_{\R} f(y)J_{\epsilon}(t,y) dy =   \int_{\mathbb{R}} f(y) L_t^ydy, \quad f\in\mathcal C^0(\R).
\end{equation*}
Ce qui  suggère de montrer la convergence de $J_{\epsilon}(t,y)$ vers $L_t^y$, quand $\epsilon$ tend vers 0. Compte tenu de (\ref{definitionJ}), cette question est équivalente à   $[ \indi_{\{ y<X_{.}\}}, X_.]_t = L_t^y$.  \\
Pour simplifier les notations, on prend $y=0$ et on note simplement $J_{\epsilon}(t)=J_{\epsilon}(t,0)$.

\section{Convergence de $J_{\epsilon}(t)$}\label{sec3}

\noindent
On peut décomposer d'une manière naturelle $J_{\epsilon}(t) $ en une somme de deux termes:
\begin{equation}
 \label{decomp1}
J_{\epsilon}(t) =
- \underbrace{ \frac{1}{\epsilon}\int_0^t
 \indi_{\{ 0<X_{s}\}}  \left( X_{s+\epsilon}-X_{s}  \right)ds}_{I^1_\epsilon(t)} +
\underbrace{ \frac{1}{\epsilon}\int_0^t \indi_{\{ 0<X_{s+\epsilon}\}}
 \left( X_{s+\epsilon}-X_{s}  \right)ds}_{I^2_\epsilon(t)}.
\end{equation}
\begin{theoreme}
\label{sec3theo1} On suppose que $X$ est un mouvement brownien
standard réel. Alors :
\begin{enumerate}
  \item $\lim_{\epsilon \to 0} \mathrm{ (ucp) }\: J_{\epsilon}(t) = [ \indi_{\{ 0<X_{s}\}}, X_s]_t = L_t^0.$
  \item $ \lim_{\epsilon \to 0} \mathrm{ (ucp) }\: I^1_{\epsilon}(t)=\int_0^t \indi_{\{ 0<X_{s}\}} dX_s$.
  \item $\lim_{\epsilon \to 0} \mathrm{ (ucp) } \: I^2_{\epsilon}(t)=X_t^+ +  \frac{1}{2} L_t^0.$
\end{enumerate}
\end{theoreme}

\medskip
\noindent \textbf{Preuve du Théorème \ref{sec3theo1}.} Puisque $f: x \to \indi_{\{ 0<x\}}\in L^2_{loc}$ et $x\to x^+$ est  une primitive de $f$,  le Théorème 4.1 de \cite{5} s'applique: $[f(X),X]_t$ existe, vaut $ 2 \left( X_t^+- X_0^+ -\int_0^t \indi_{\{ 0<X_{s}\}} dX_s\right)$ et la formule de Tanaka donne alors le point $(i)$. Toujours par  \cite{5}, on a
\begin{equation}
\label{sec3label1}
\lim_{\epsilon \to 0}\; \textrm{ (ucp) } \frac{1}{\epsilon}\int_0^t
 Y_s  \left( X_{s+\epsilon}-X_{s}  \right)ds=\int_0^t
 Y_s  dX_s,
\end{equation}
avec $Y_s=f(X_s)$ , ce qui donne le point $(ii)$. Le point $(iii)$ se déduit des deux précédents via la formule de Tanaka. 
\qed

\medskip
\noindent \textbf{Remarque:}  
Plus généralement, si $(X_t)_{t\geqslant 0}$ est une semi-martingale et  $(Y_t)$ un processus adapté tel que $t\to Y_t$ admet des limites à gauche, alors  (\ref{sec3label1}) a lieu (c.f. \cite{8}, Proposition 3.31). Signalons un résultat qui ne concerne pas directement  l'approximation du temps local mais qui est très lié à notre étude: si $(Y_t)$ est un processus adapté et localement h\"oldérien, alors la convergence  (\ref{sec3label1}) a lieu presque s\^urement, uniformément pour $t\in [0,T]$.

\section{Autres schémas d'approximation}\label{sec4}

\noindent D'après la décomposition (\ref{decomp1}) et le Théorème \ref{sec3theo1}, $J_{\epsilon}(t)$ se décompose en deux termes ayant chacun une limite. Mais ces deux limites ne s'expriment pas uniquement en fonction du temps local. On cherche donc une autre  décomposition de $J_{\epsilon}(t)$ en des termes qui convergent chacun vers une fraction du temps local. En écrivant $ \indi_{\{X_{(u+\epsilon)\wedge t } >0\}} -\indi_{\{X_{u} >0\} }= \indi_{\{X_{(u+\epsilon)\wedge t  } >0,   X_u \leqslant 0\} }-\indi_{\{X_{(u+\epsilon)\wedge t } \leqslant 0,   X_u > 0\} }$, on obtient facilement:
\begin{equation}
\label{decomp2}
J_{\epsilon}(t) =I^3_\epsilon(t) + I^4_\epsilon(t) + I^5_\epsilon(t) + R_\epsilon(t),
\end{equation}
\begin{eqnarray}
\textrm{avec }\qquad I^3_\epsilon(t) & =&
 \frac{1}{\epsilon}\int_0^{t} X_{(u+\epsilon)\wedge t}^+ \indi_{ \{X_u \leqslant 0\} } du
+
\frac{1}{\epsilon}\int_0^{t} X_{(u+\epsilon) \wedge t}^- \indi_{ \{X_u>0\} } du, \label{decomp2b}\\
 I^4_\epsilon(t)&=& \frac{1}{\epsilon}\int_0^{t} X_u^- \indi_{ \{X_{(u+\epsilon) \wedge  t} > 0\} } du, \qquad
  I^5_\epsilon(t)  =\frac{1}{\epsilon}\int_0^{t} X_u^+ \indi_{ \{X_{(u+\epsilon) \wedge t} \leqslant 0\} } du,\label{decomp2d} 
\end{eqnarray}
et $R_{\epsilon}(t)$ un terme qui converge presque sûrement vers 0, uniformément sur les compacts. 

   \medskip
\begin{theoreme}
\label{sec4eq1}
Si $X$ est le mouvement Brownien standard, alors
\begin{enumerate}
  \item $ \lim_{\epsilon \to 0} \mathrm{ (ucp) } \:  I^3_{\epsilon}(t)= \frac{1}{2} L_t^0.$
  \item $ \lim_{\epsilon \to 0} \mathrm{ (ucp) } \:  I^4_{\epsilon}(t)=  \lim_{\epsilon \to 0}\: \mathrm{ (ucp) } \:  I^5_{\epsilon}(t)= \frac{1}{4} L_t^0.$
\end{enumerate}
\end{theoreme}

\medskip
\noindent \textbf{Remarque:}  
\begin{itemize}
\item[1)] Ce résultat est encore vrai lorsque $X_t=\int_0^t \sigma(s)dB_s$ où $B$ est un mouvement brownien standard et $\sigma$ une fonction  définie sur $\R^+$, h\"olderienne d'ordre $\gamma > \frac{1}{4}$ et  telle que $| \sigma(s) | \geqslant  a > 0$ .
\item[2)] Nous n'avons pas obtenu séparément la convergence de chacun des termes de $I^3_{\epsilon}(t)$.
\end{itemize}

\medskip
\noindent \textbf{Preuve du Théorème \ref{sec4eq1}.} \textbf{a)} D'après la formule de Tanaka,
$$X_{(u+\epsilon)\wedge t }^+= X_u^+ +\int_u^{(u+\epsilon)\wedge t }\indi_{\{X_s > 0\} } dX_s + \frac{1}{2}(L_{(u+\epsilon)\wedge t }^0-L_u^0).$$
On a une formule similaire pour $X_{(u+\epsilon)\wedge t }^-$. Il est aisé d'en déduire que $I^3_{\epsilon}(t)$ est égal à
\begin{equation}
\label{sec3eq4}
 \frac{1}{\epsilon}\int_0^{t} \left( \int_u^{(u+\epsilon)\wedge t } \left(\indi_{\{X_s > 0, X_u \leqslant 0\} }-\indi_{\{X_s \leqslant 0, X_u > 0\}}\right) dX_s \right) du  +\frac{1}{2\epsilon}\int_0^{t} (L_{(u+\epsilon)\wedge t }^0-L_u^0)du.
\end{equation}
On écrit
$L_{(u+\epsilon)\wedge t }^0-L_u^0= \int^{(u+\epsilon)\wedge t }_u dL^0_s$.  Une application
 du théorème de Fubini permet montrer la convergence du second terme de (\ref{sec3eq4}) vers $\frac{1}{2}L^0_t$. On  utilise  le théorème de Fubini stochastique (c.f. Section IV.5 de \cite{1}) pour transformer le premier terme de (\ref{sec3eq4}) en  une intégrale stochastique. L'inégalité de Doob donne alors la convergence  de ce terme vers 0 dans $L^2(\Omega)$.

\noindent \textbf{b)} Il est équivalent d'étudier la convergence de $
I^4_{\epsilon}(t)$ ou  de $ I^5_{\epsilon}(t)$. On \'ecrit $ I^5_{\epsilon}(t)$ comme la somme d'un terme
 qui converge presque  sûrement vers 0, uniformément sur les compacts, et de:
\begin{equation}
\label{sec3eq3}
 \frac{1}{\epsilon} \int_{0}^{(t-\epsilon)^+} X_{u}^+
 \left[ \indi_{\{ X_{u+\epsilon}\leqslant 0\}}-
 \Phi \left(- \frac{X_u}{\sqrt{\epsilon}}\right)\right] du +
 \frac{1}{\epsilon} \int_{0}^{t} X_{u}^+ \Phi \left(-
 \frac{X_u}{\sqrt{\epsilon}}\right)du,
\end{equation}
o\`u  $\Phi$ est la fonction de répartition de la gaussienne
centrée réduite. Par la formule de densité d'occupation:
 $$  \frac{1}{\epsilon} \int_{0}^{t} X_{u}^+ \Phi \left(-
 \frac{X_u}{\sqrt{\epsilon}}\right)du = \frac{1}{\epsilon} \int_{\R}  x^+ \Phi \left( -\frac{x}{\sqrt{\epsilon}}\right)L^x_t dx = \int_0^{\infty}   y \Phi ( -y  ) L^{y\sqrt{\epsilon}}_t    dy.$$
On en déduit facilement la convergence p.s. de ce terme vers  $\frac{1}{4} L_t^0$, uniformément sur les compacts.  Pour le premier terme de (\ref{sec3eq3}),  on écrit le terme entre crochet comme une intégrale stochastique, puis on utilise le théorème de Fubini pour se ramener une martingale.  Une utilisation de l'inégalité de Doob permet d'obtenir la convergence vers 0 dans $L^2(\Omega)$, uniformément sur les compacts.\qed

\section{Vitesse de convergence de $J_{\epsilon}(t)$ dans $L^2(\Omega)$} \label{sec5}

\begin{theoreme}
\label{sec5eq1}
Soit $(X_t)$  le mouvement brownien standard. Pour tout $T>0$, $\delta \in ]0,\frac{1}{4}[$, il existe une constante $K_{\delta}$ telle que:
\begin{equation}
\label{sec5eq1b}
 \forall \epsilon \in ]0,1], \quad \left\| \sup_{t\in [0,T]} \left| J_{\epsilon}(t) -  L_t^0 \right|  \right\|_{L^2(\Omega)}  \leqslant K_{\delta}\epsilon^{\delta}.
\end{equation}
\end{theoreme}
On a un résultat similaire pour la vitesse de convergence de  $I^i_{\epsilon}(t)$ vers sa limite, $i=1,\dots, 5$.

   \medskip
\noindent \textbf{Preuve du Théorème \ref{sec5eq1}.} On utilise des éléments des preuves précédentes:
\begin{equation}
\label{sec5eq2}
J_{\epsilon}(t) -  L_t^0= -\left(I^1_{\epsilon}(t)-  \int_0^t \indi_{\{ 0<X_{s}\}} dX_s \right) +  \left( I^4_{\epsilon}(t) -\frac{1}{4}L_t^0(X) \right)  + \left( I^5_{\epsilon}(t) -\frac{1}{4}L_t^0(X) \right) +R_{\epsilon}(t) \end{equation}
On décompose chaque terme sous la forme $\int_0^t h(s,\epsilon) dX_s + \int_0^t k(s,\epsilon) ds$.  En utilisant la propriété de Hölder du mouvement Brownien et de son temps local, on majore p.s. $\sup_{t\in [0,T]}  |\int_0^t  k(s,\epsilon) ds|$ par $C\epsilon^\delta$. Grâce à l'inégalité de Doob, on majore $E(( \sup_{t\in [0,T]}  \int_0^t h(s,\epsilon) dX_s)^2)$ par $4\int_0^T E(  h^2(s,\epsilon)) ds$. Il est possible, après des calculs longs et fastidieux, de montrer que ce terme est lui-même majoré par  $C\sqrt{\epsilon}$.
\qed

\medskip
\noindent \textbf{Remarque:}  
Malheureusement, ($\ref{sec5eq1b}$) ne permet pas de montrer la convergence p.s. Il est toutefois possible, en modifiant la preuve du Théorème \ref{sec5eq1} et en utilisant le lemme de  Borel-Cantelli, de montrer que $\sup_{t\in [0,T]} \left| J_{\epsilon_n}(t) -  L_t^0 \right| $  converge presque sûrement vers 0, lorsque $n\to \infty$, où $ (\epsilon_n)_{n\in \N}$ est une suite  réelle positive décroissante tendant vers 0 et telle que  $\sum_{i=1}^{\infty} \sqrt{\epsilon_i}<\infty$.

{\small
\def\refname{R\'ef\'erences}
\bibliography{refCRASBBPV}
\bibliographystyle{plain}
}

\newpage

\noindent Dear Editor,

\medskip
\noindent First we would like to gratefully acknowledge  the referee for his (her) suggestions.

\medskip
 \noindent According to referee's remark, we have changed the abstract as follows:
 
\medskip
 \begin{itemize}  
\item[] 
{\bf Some Brownian local time approximations} \\
 We give some approximations of  the local time process $(L_t^x)_{t\geqslant 0}$ at level $x$ of the real Brownian motion $(X_t)$.
 We prove  that $ \frac{2}{\epsilon}\int_0^{t} X_{(u+\epsilon)\wedge t}^+ \indi_{ \{X_u \leqslant 0\} } du +
\frac{2}{\epsilon}\int_0^{t} X_{(u+\epsilon) \wedge t}^- \indi_{ \{X_u>0\} } du$ 
and $\frac{4}{\epsilon}\int_0^{t} X_u^- \indi_{ \{X_{(u+\epsilon) \wedge  t} > 0\} } du$ converge  in the ucp sense to $L_t^0$, as $\epsilon \to 0$. 
 We show  that $ \frac{1}{\epsilon}\int_0^t ( \indi_{\{ x<X_{s+\epsilon}\}} - \indi_{\{ x<X_{s}\}} )   ( X_{s+\epsilon}-X_{s} )ds$ 
 goes  to $L_t^x$  in $L^2(\Omega)$   as $\epsilon \to 0$, and that the rate of convergence is of order $\epsilon^\alpha$, for any $\alpha < \frac{1}{4}$. 
\end{itemize}

\medskip
\noindent  The title of Section 2. has been changed.

\noindent The mistakes on the world "h\"olderien" has been corrected

\noindent  The bibliography has been homogenized to follow the alphabetic order.

\vskip 10pt
\noindent Yours sincerely.

\vskip 10pt \noindent Nancy on  20th April 2007,

 \vskip 10pt \noindent  B. B\'erard-Bergery, \quad P. Vallois.

\end{document}